\documentclass[12pt]{article}
\usepackage{graphicx,amssymb,amsmath}
\usepackage{hyperref,url}
\usepackage{braket}
\usepackage{verbatim}
\newcommand{\nocopyright}{
No Copyright\thanks{
The authors hereby waive all copyright
and related or neighboring rights to this work,
and dedicate it to the public domain.
This applies worldwide.
}}
\title{A category for bijective combinatorics}
\author{Peter G. Doyle}
\date{Version 16 dated 20 June 2019\\
\nocopyright
}
\newtheorem{theorem}{Theorem}

\newtheorem{corollary}[theorem]{Corollary}

\newcommand{\proofstart}{{\noindent \bf Proof.\ }}
\newcommand{\qed}{\spadesuit}
\newcommand{\mathproofend}{\quad \qed}
\newcommand{\proofend}{$\quad \qed$}

\newcommand{\p}{\;\;\;}
\newcommand{\comma}{\p,}
\newcommand{\pd}{\p.}

\newcommand{\fig}[3]{
\begin{figure}
\includegraphics[width=370pt]{figures/#1.pdf}
\caption{#3}
\label{#2}
\end{figure}
}
\newcommand{\sub}[1]{\braket{#1}}
\newcommand{\nnn}{\mathbf{NSet}}
\newcommand{\zzz}{\mathbf{ZSet}}
\newcommand{\nn}{\mathbf{N}}
\newcommand{\zz}{\mathbf{Z}}

\makeatletter
\newcommand{\xRrightarrow}[2][]{\ext@arrow 0359\Rrightarrowfill@{#1}{#2}}
\newcommand{\Rrightarrowfill@}{\arrowfill@\equiv\equiv\Rrightarrow}
\newcommand{\xLleftarrow}[2][]{\ext@arrow 3095\Lleftarrowfill@{#1}{#2}}
\newcommand{\Lleftarrowfill@}{\arrowfill@\Lleftarrow\equiv\equiv}
\newcommand{\xLleftRrightarrow}[2][]{\ext@arrow 3399\LleftRrightarrowfill@{#1}{#2}}
\newcommand{\LleftRrightarrowfill@}{\arrowfill@\Lleftarrow\equiv\Rrightarrow}
\makeatother

\newcommand{\eq}{\Rrightarrow}
\newcommand{\funeq}[1]{\stackrel{#1}\eq}

\newcommand{\union}{\cup}
\newcommand{\cross}{\times}
\newcommand{\op}{\operatorname}
\newcommand{\kw}{\operatorname}

\newcommand{\cancel}{\op{cancel}}
\newcommand{\nestuntil}{\operatorname{nestuntil}}

\newcommand{\id}{\operatorname{id}}
\newcommand{\create}{\operatorname{create}}
\newcommand{\destroy}{\operatorname{destroy}}
\newcommand{\then}{\triangleleft}

\begin{document}

\maketitle

\begin{abstract}
The category of matchings between finite sets
extends to the category of cobordisms of signed sets.
A chain of cobordisms that starts and ends with
unsigned sets $A$ and $B$
yields a matching from $A$ to $B$.
This is a convenient way to package the involution principle of
Garsia and Milne,
which reveals itself to have little to do with involutions.
\end{abstract}

\section{Introduction}

As observed in passing by Conway and Doyle
\cite[p.\ 23]{conwaydoyle:three},
and doubtless by others before them,
bijective combinatorics can be viewed as cobordism theory for
oriented 0-dimensional manifolds.
We develop this approach,
with a view to clarifying the role of the involution principle
of Garsia and Milne
\cite{garsiaMilne:method,garsiaMilne:rogers}.

There will be nothing new here,
beyond notation.
Cobordism theory dates from the 1950s,
but its 0-dimensional manifestations can be seen in
what is now the standard
proof of the Cantor-Schroeder-Bernstein equivalence theorem,
given by Koenig \cite{koenig} in 1906.
(See Appendix \ref{app}.)
And some will see the origins
even further back in the mists of time.
The application to combinatorics is implicit in
Picciotto
\cite{picciotto:thesis},
and hardly different from the approach taken in texts like
Stanton and White
\cite{stantonwhite},
and indeed the papers
Garsia and Milne.
In the end, it all comes down to subtraction.

\section{Notation}

We deal with \emph{matchings} (bijections) between finite sets.
To emphasize this, we write
\[
f: A \eq B
\]
(`$f$ matches A to B'), or
\[
A \funeq{f} B
\pd
\]

We write
composition in natural order,
using the symbol $\then$, pronounced `then':
\[
(f \then g)(x) = g(f(x))
\pd
\]

\begin{comment}
In place of a function name we may put a recipe (say, a $\lambda$-expression),
so that for example we could
write
\[
\lambda x.4x: \zz/17 \eq \zz/17
\comma
\]
and say, `Multiplication by four permutes the integers mod seventeen.'.
\end{comment}

We write $X + Y$ for the disjoint union
\[
X+Y = X \cross \{0\} \union Y \cross \{1\}
\p,
\]
and adopt all the usual type coercions (`abuses of notation'),
so that
\[
X \subset X + Y
;\;
X+Y = Y+X
;\;
(X+Y)+Z = X+(Y+Z)
\p,
\]
etc.
We can also take the disjoint union of matchings:
If
\[
f: A \eq B
\]
and
\[
g: C \eq D
\]
then
\[
f+g: A+B \eq C+D
\pd
\]

\section{Simple subtraction}

\begin{theorem}[Simple subtraction] \label{simplesubtraction}
If
\[
f: A + C \eq B + C
\]
then
\[
\cancel(C,f):
A \eq B
\comma
\]
where
\[
\cancel(C,f)(a) =
\nestuntil(\lambda x.x \notin C,f)
(f(a))
\]
and
\[
\op{nestuntil}(\op{test},f) = 
\lambda x.
\kw{if}
\op{test}(x)
\kw{then}
x
\kw{else}
\nestuntil(\op{test},f)
(f(x))
\pd
\]
\end{theorem}

\proofstart
Figure \ref{fig:simplesubtraction} shows the idea.
The pidgin $\lambda$-calculus used to define $\cancel(C,f)$ just means that,
given $a$,
we start with $f(a)$ and keep applying $f$ until the result
escapes from $C$.
This happens eventually because $C$ is finite.
This escape mechanism applies to any function $f:A+C \to B+C$
to yield a function $\cancel(C,f):A \to B$.
If $f$ is an injection, so is $\cancel(C,f)$;
if $f$ is a surjection, so is $\cancel(C,f)$.
\proofend
\fig{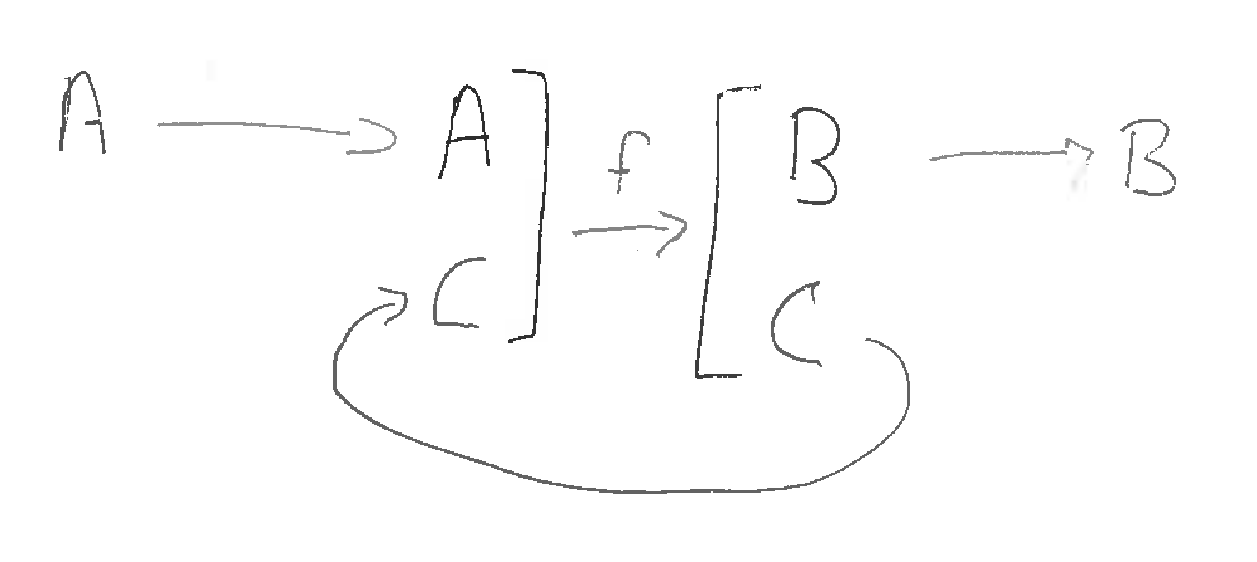}{fig:simplesubtraction}{Simple subtraction.}

\vspace{0.25cm}
\noindent
{\bf Aside.} We don't need $A,B,C$ to be finite here:
$A$ and $B$ may be infinite,
and $C$ need only be `Dedekind-finite'---but
let's not open that
can of worms.

\begin{corollary}[Repeated subtraction]\label{repeatedsubtraction}
Let
\[
f: A+(C \union D) \to B+(C \union D)
\]
where $C,D$ are finite sets, not necessarily disjoint.
Then
\[
\cancel(D,\cancel(C,f))
=
\cancel(C \union D,f)
=
\cancel(C,\cancel(D,f))
:A \eq B
\pd
\]
\end{corollary}

\proofstart
If you persist in escaping from $C$ until you have left $D$,
to an outside observer your behaviour is indistinguishable from
escaping from $C \union D$.
\proofend

\vspace{0.5cm}

Here ends the math.
The rest is bookkeeping.

\newpage
\section{The cobordism category}

We've been working in the category
whose morphisms are matchings of finite sets.
Let's call this category $\nnn$.
Paralleling the extension of the natural numbers $\nn$ to the integers $\zz$,
we are going the extend $\nnn$ to the category $\zzz$ whose morphisms
are cobordisms of signed sets.

\subsection{Signed sets}

Define a \emph{signed set} to be an ordered pair 
$A = \sub{A^+,A^-}$ of finite sets.
We identify unsigned sets as signed sets $A$ for which $A^-=0$,
where $0$ is the empty set, so that
$\sub{A^+,0}=A^+$,
and in particular $\sub{0,0} = 0$.
Soon we'll be writing $A = A^+ - A^-$,
which is how we want to think of it.

Write
\[
|A|=A^+ + A^-
\comma
\]
and define
\[
A \subset B \iff
|A| \subset |B|
\iff
A^+ \subset B^+ \land A^- \subset B^-
\pd
\]
If $A \subset B$, define
\[
B \setminus A = \sub{B^+ \setminus A^+,B^- \setminus A^-}
\subset B
\pd
\]

Define
\[
A \leq B
\iff
A^+ \subset B^+ \land A^- \supset B^-
\comma
\]
and observe that being unsigned is the same as being (weakly) positive:
\[
A \geq 0
\iff
\mbox{$A$ is unsigned}
\iff
A = A^+
\iff
A = |A|
\pd
\]
(This way of designating unsigned sets is the only use we will make
of this partial order.)

Define sum, negation, and difference of signed sets in the obvious ways:
\[
A+B = \sub{A^+ + B^+, A^- + B^-}
\p;
\]
\[
-A = \sub{A^-,A^+}
\p;
\]
\[
A-B = A+(-B) = \sub{A^+ + B^-, A^- + B^+}
\pd
\]
Now for a signed set $A$ we can write
\[
A = \sub{A^+,A^-} = A^+ - A^-
\pd
\]
Everything works as expected, except
that we can't replace $A-A$ with $0$, or vice versa.

\subsection{Cobordisms}

A \emph{cobordism} is a matching of signed sets,
represented by a triple
$(|f|,A,B)$, where $A,B$ are signed sets, and
\[
|f|: A^+ + B^- \eq A^- + B^+
\pd
\]
We use the same notation for
cobordisms as for matchings, writing
\[
f: A \eq B
\]
and saying `f matches $A$ to $B$' (or `f is a cobordism from $A$ to $B$',
if we want to emphasize that we're dealing with signed sets).

If
\[
A \funeq{f} B \funeq{g} C
\comma
\]
we have matchings
\[
|f|: A^+ + B^- \eq A^- + B^+
\comma
\]
\[
|g|: B^+ + C^- \eq B^- + C^+
\pd
\]
These combine to give a matching
\[
|f|+|g|
:
A^+ + |B| + C^- \eq A^- + |B| + C^+
\pd
\]
We define
\[
f \then g: A \eq C
\]
by setting
\[
|f \then g| = \cancel(|B|,|f|+|g|): A^+ + C^- \eq A^- + C^+
\pd
\]
(Cf.\ The `Bread Lemma' of Picciotto
\cite[p.\ 25, Lemma 2]{picciotto:thesis}).

From Corollary \ref{repeatedsubtraction} (repeated subtraction) we get
\begin{corollary}[Associativity for composition of cobordisms]
Suppose
\[
A \funeq{f} B \funeq{g} C \funeq{h} D.
\]
Then
\[
((f \then g) \then h) = (f \then (g \then h)) = f \then g \then h:
A \eq D
\comma
\]
where
\[
|f \then g \then h| = \cancel(|B|+|C|,|f|+|g|+|h|): A^+ + D^- \eq A^- + D^+
\pd
\mathproofend
\]
\end{corollary}

\begin{corollary}[Chain associativity]
If
\[
A_0 \funeq{f_{0,1}}
A_1 \funeq{f_{1,2}}
\ldots \funeq{f_{n-2,n-1}}
A_{n-1} \funeq{f_{n-1,n}} A_n
\]
then
\[
|f_{0,1} \then \ldots \then f_{n-1,n}|
=
\cancel(|A_1|+\ldots+|A_{n-1}|,|f_{0,1}|+\ldots+|f_{n-1,n}|)
\pd
\]
In particular, if $A_0,A_n \geq 0$ then
\[
\cancel(|A_1|+\ldots+|A_{n-1}|,|f_{0,1}|+\ldots+|f_{n-1,n}|)
:A_0 \eq A_n
\mathproofend
\]
\end{corollary}
(Cf. Picciotto
\cite[p.\ 26, Lemma 3]{picciotto:thesis}.)
We've written this out to emphasize that composing a chain of cobordisms 
requires only a single application of subtraction.

Of course we also have identity cobordisms
\[
\id(A):A \eq A
\]
with
\[
|\id(A)|:=\id(|A|)
\comma
\]
where in the second instance $\id$ denotes the identity in $\nnn$.
So we have ourselves a category, which we call $\zzz$.

Just as unsigned sets correspond to signed sets $A$ with $A \geq 0$,
matchings correspond to cobordisms
\[
f: A \eq B
\]
where $A,B \geq 0$, so that
\[
|f|: A \eq B
\pd
\]
This correspondence is natural (`functorial'):
If $A,B,C \geq 0$ and
\[
A \funeq{f} B \funeq{g} C
\]
then
\[
|f \then g| = |f| \then |g|: A \eq C
\pd
\]
So we can identify a matching as a cobordism $f$ for which $f=|f|$,
just as an unsigned set is a signed set for which $A = |A|$.

\subsection{Arithmetic with cobordisms}

We can add, negate, and subtract cobordisms, 
bearing in mind that negating a cobordism reverses the direction of the arrow:
If
\[
f:A \eq B
\]
then
\[
-f:-B \eq -A
\comma
\]
with
\[
|-f|=|f|:B^- + A^+ \eq B^+ + A^-
\pd
\]

Now we can write the identity cobordism as
\[
\id(A) = \id(A^+) - \id(A^-): A \eq A
\pd
\]
Closely related to the identity are the creation and destruction cobordisms
\[
\create(A): 0 \eq A-A
\]
and
\[
\destroy(A) = -\create(A): A-A \eq 0
\comma
\]
where
\[
|\create(A)| = |\destroy(A)| = |\id(A)| = \id(|A|)
\pd
\]
Observe that
\[
\create(A) = \create(-A) = \create(|A|): 0 \eq A-A = (-A) - (-A) = |A| - |A|
\pd
\]

More generally, for any
\[
\phi: A \eq B
\]
we define
corresponding creation and destruction morphisms
\[
\create(\phi):0 \eq B-A
\]
and
\[
\destroy(\phi) = -\create(\phi): A - B \eq 0
\]
with
\[
|\create(\phi)| = |\destroy(\phi)| = |\phi|: A^+ + B^- \eq A^- + B^+
\pd
\]
We have
\[
\create(\phi) = \create(-\phi) = \create(|\phi|): 0 \eq B-A
\pd
\]

\section{Signed subtraction} \label{sec:signedsub}

With this machinery in place,
we immediately get
\begin{corollary}[Subtraction] \label{subtraction}
For any signed sets $A,B,C,D$, if
\[
f:A+C \eq B+D
\]
and
\[
g:D \eq C
\]
then
\[
(\id(A) + \create(C)) \then (f-g) \then (\id(B) - \create(D)): A \eq B
\pd
\]
\end{corollary}

\proofstart
\[
-g:-C \eq -D
\comma
\]
so
\[
f-g:A+C-C \eq B+D-D
\comma
\]
so
\[
A
\funeq{\id(A)+\create(C)}
A+C-C
\funeq{f-g}
B+D-D
\funeq{\id(B)-\create(D)}
B
\pd
\mathproofend
\]
Figure \ref{fig:subtraction} shows a diagram.
\fig{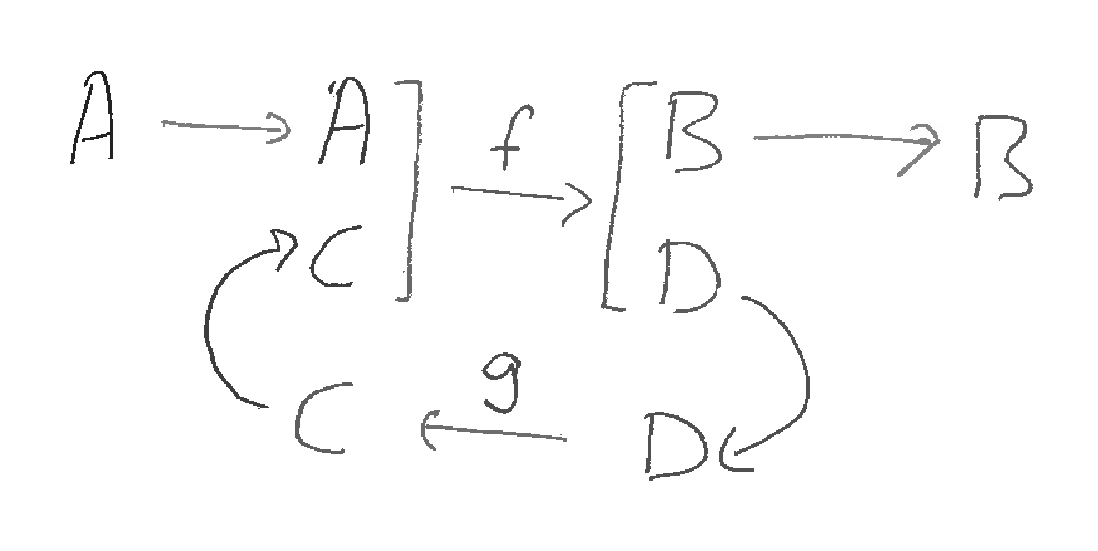}{fig:subtraction}{Subtraction.}
Taking $C=D$, $g=\id(C)$ we recover Figure \ref{fig:simplesubtraction},
so it might seem that we've made little progress beyond
Theorem \ref{simplesubtraction},
and in a sense this is very true.
Of course we now have subtraction working for signed sets,
as emphasized in the exploded view of Figure \ref{fig:signedsub}.
\fig{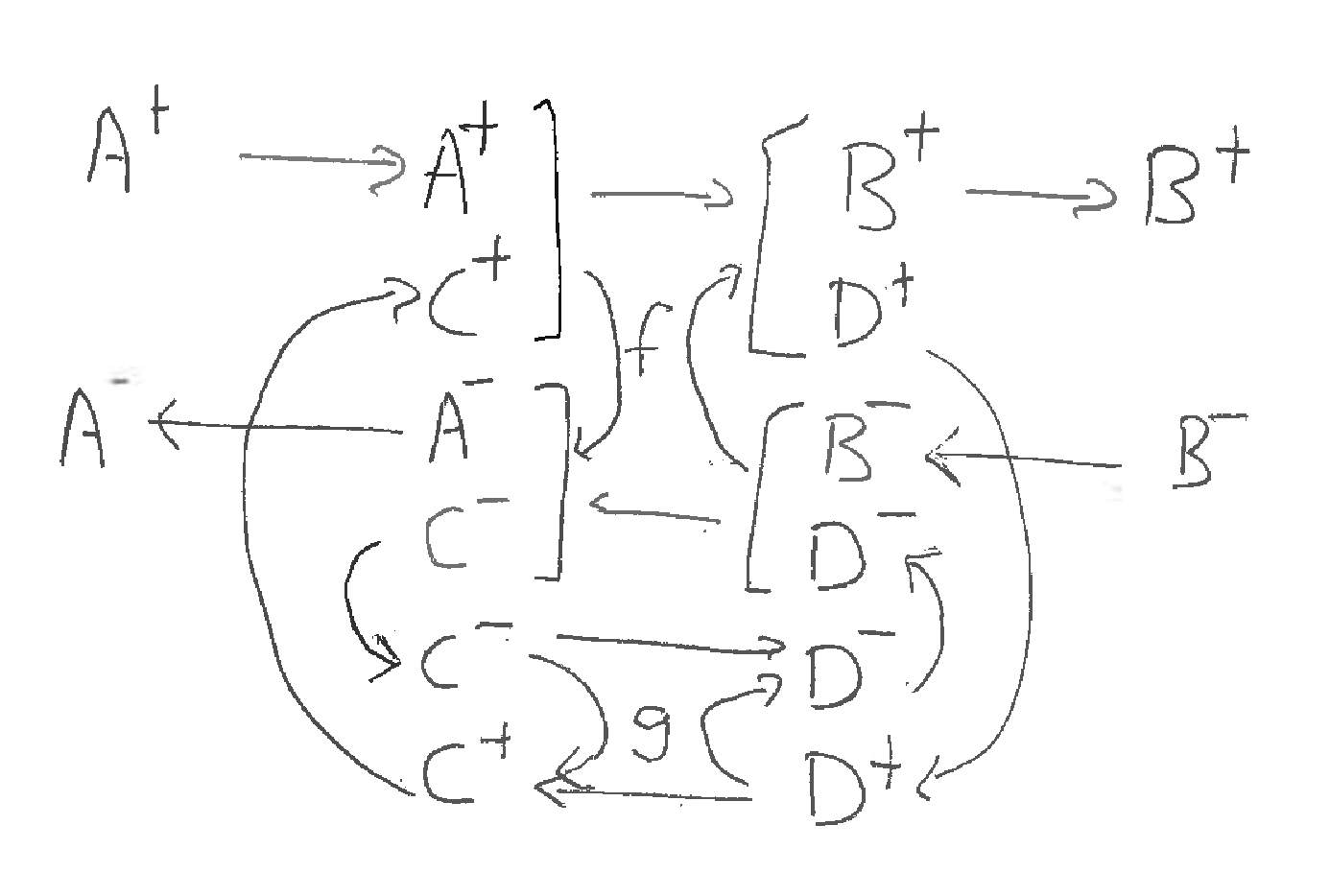}{fig:signedsub}{Subtraction (exploded view).}
But mainly, we've just taken simple subtraction,
dressed it up, and called it `composition of morphisms'.

\section{The involution principle}
\begin{corollary}[The involution principle]
For any signed sets $X,Y$ and $A,B \subset X$, if
\[
\phi:Y \eq X \setminus A
\]
and
\[
\psi:
X \setminus B \eq Y
\]
then
\[
(\id(A) + \create(\phi)) \then (\id(B) - \create(\psi)): A \eq B
\pd
\]
\end{corollary}

\proofstart
\[
\create(\phi): 0 \eq (X \setminus A) - Y
\]
and
\[
-\create(\psi): (X \setminus B) - Y \eq 0
\]
so
\[
A
\funeq{\id(A) + \create(\phi)}
A+ (X \setminus A) - Y
=
X-Y
=
B + (X \setminus B) - Y
\funeq{\id(B) -\create(\psi)}
B
\pd
\mathproofend
\]
Figure \ref{fig:involution} shows the diagram;
\fig{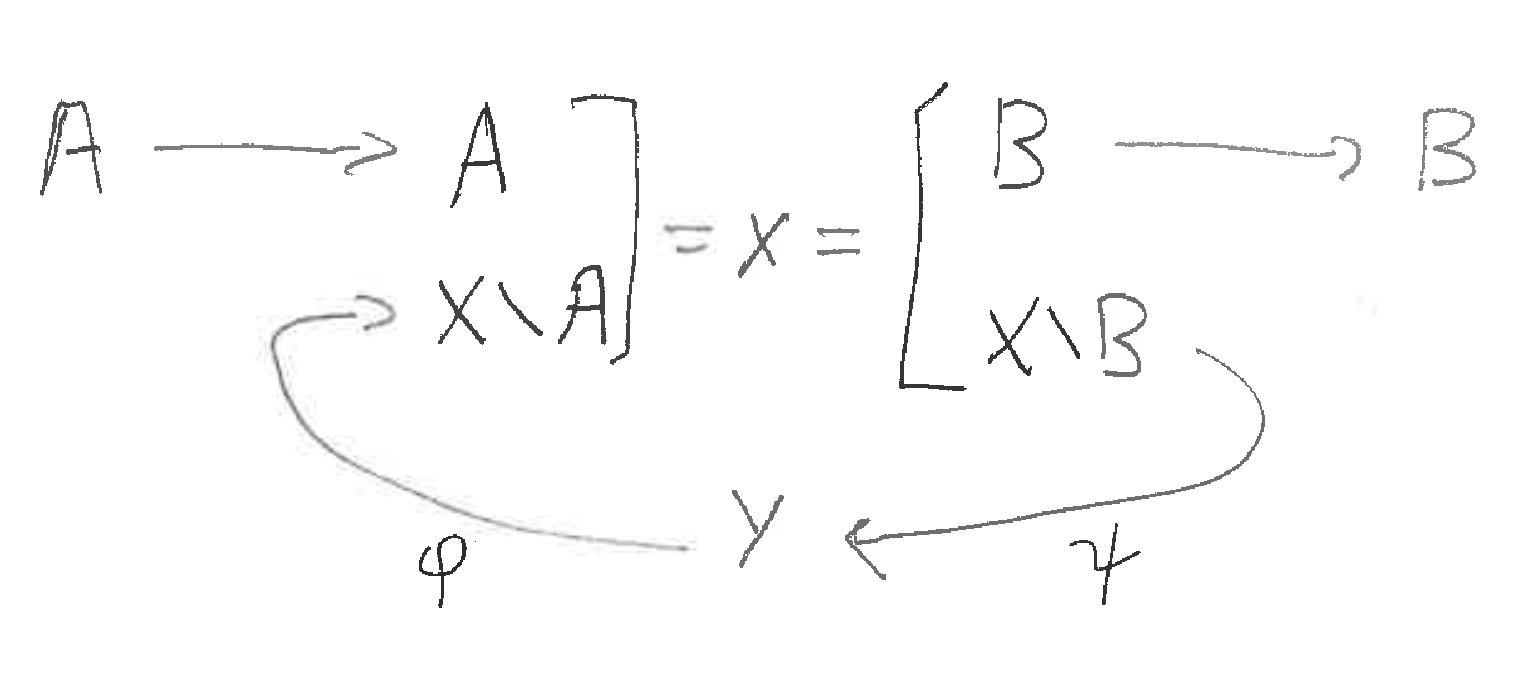}{fig:involution}{The involution principle.}
Figure \ref{fig:involutionexploded} shows the exploded view.
\fig{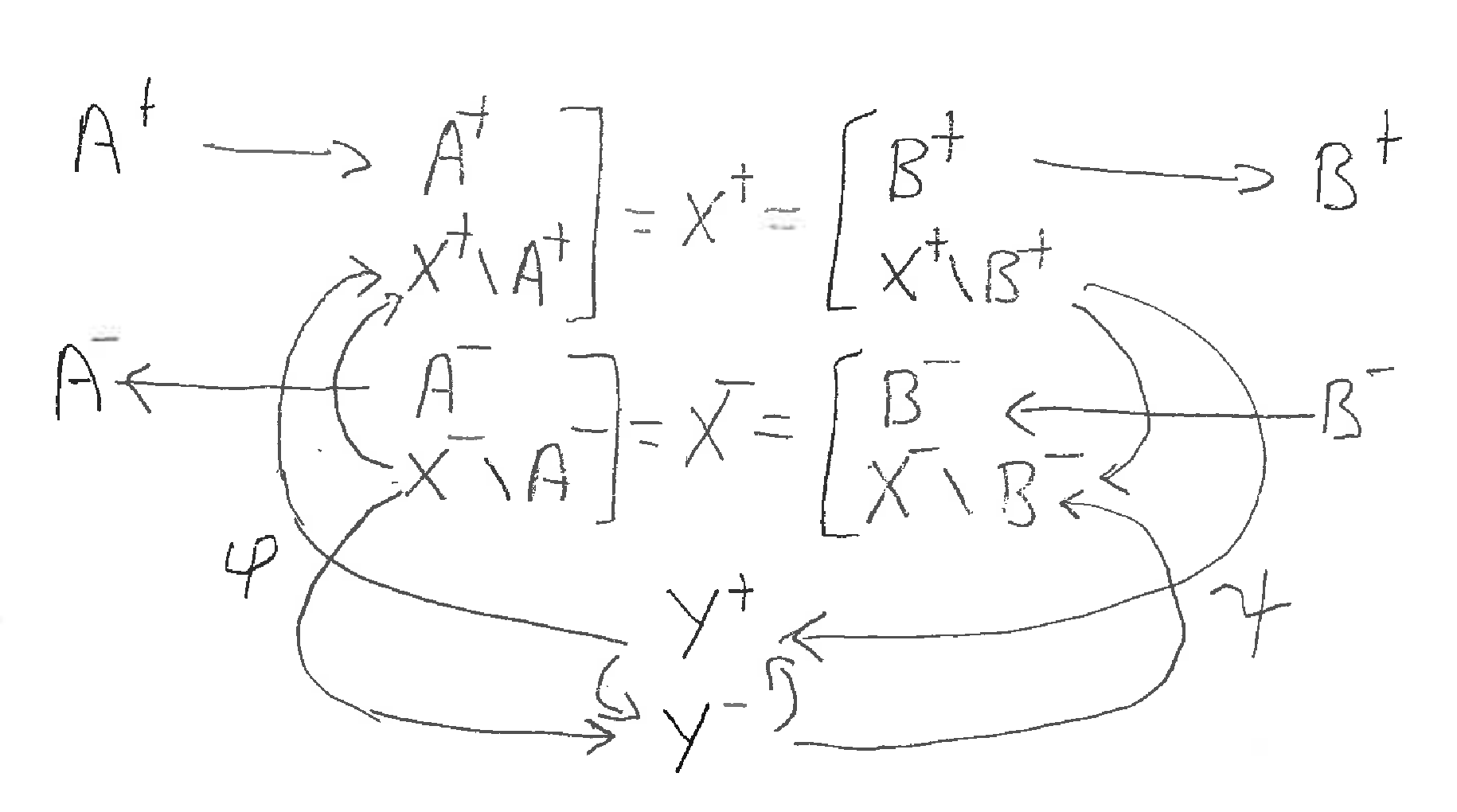}{fig:involutionexploded}{The involution principle (exploded view).}

If we restrict to unsigned sets $X,Y \geq 0$ we recover
the Garsia-Milne involution principle.
\begin{corollary}[The Garsia-Milne involution principle]
If $A \subset X \geq 0$, $B \subset Y \geq 0$,
$\phi:Y \eq X \setminus A$,
$\psi:X \setminus B \eq Y$
then
\[
h =(\id(A) + \create(\phi)) \then (\id(B) - \create(\psi)): A \eq B
\comma
\]
where
\[
h(a) =
\nestuntil(\lambda x.x \in B,\lambda x.\phi(\psi(x)))
(a)
\pd
\mathproofend
\]
\end{corollary}

Formulated in this way, the involution principle has no need of involutions.
If we want them nevertheless,
we can manufacture them,
as long as we're willing to take inverses of matchings.
(Have you noticed that we've been studiously avoiding doing this?)
From
\[
\phi: Y \eq X \setminus A
\]
we get
\[
\phi^{-1}: X \setminus A \eq Y
\comma
\]
so
\[
\Phi = \id(A) + \phi + \phi^{-1}:
X+Y = A + Y + (X \setminus A)
\eq
A + (X \setminus A) + Y
=
X+Y
\comma
\]
with
\[
\Phi \then \Phi = \id(X+Y)
\pd
\]
Likewise,
\[
\psi: X \setminus B \eq Y
\comma
\]
\[
\psi^{-1}: Y \eq X \setminus B
\comma
\]
\[
\Psi = \id(B) + \psi + \psi^{-1}:
X+Y = B + (X \setminus B) + Y \eq B + Y + (X \setminus B)
=
X+Y
\comma
\]
with
\[
\Psi \then \Psi = \id(X+Y)
\pd
\]
Now we have
\[
h(a) =
\nestuntil(\lambda x . \Psi(x) = x, \lambda x. \Phi(\Psi(x)))
(a)
\pd
\]

In practice, 
$\phi$ and $\psi$ will often naturally arise as restrictions of
involutions $\Phi$ and $\Psi$.
And when this happens, exchanging $\Phi$ and $\Psi$ gets us the
inverse matching:
\[
h^{-1}(b) =
\nestuntil(\lambda x . \Phi(x) = x, \lambda x. \Psi(\Phi(x)))
(b)
\pd
\]

Still, it will be helpful to recognize that fundamentally,
the involution principle has little to do with involutions.

\appendix
\section{Sidestepping division}

The paper of Garsia and Milne
\cite{garsiaMilne:rogers}
is one of the landmarks of bijective
combinatorics.
Beyond the specific application to the Rogers-Ramanujan identities,
this work was a triumph for the null hypothesis that
where there is algebra, there is combinatorics;
it introduced combinatorialists to subtraction;
and it showed the virtue of working with signed sets,
which we've been touting.

Here we call attention to yet another aspect of their work,
which was the way they avoided having to divide.

At its most basic, the problem of division is this.
Suppose $A,B,C$ are finite sets, with $C \neq 0$.
From
\[
f:A \cross C \eq B \cross C
\comma
\]
we want to produce
\[
h:A \eq B
\pd
\]
Rephrased for signed sets, 
from
\[
f: A \cross C \eq 0
\comma
\]
we want to produce either
\[
g: C \eq 0
\]
or
\[
h: A \eq 0
\pd
\]

Subtraction is straight-forward, but division is not.
There are situations where division is needed,
and techniques that will make it work.
(Cf.\ Feldman and Propp \cite{fp};
Doyle and Qiu \cite{doyleQiu:four};
Bajpai and Doyle \cite{bajpaiDoyle:equivariant}.)
Garsia and Milne's insight was that,
when working with generating functions,
multiplying by the reciprocal may obviate the need to divide.

To give the idea,
suppose we have
generating function $F,G$ that we wish to show are equal.
We have a bijection showing
$F \cdot H = G \cdot H$,
and we want to derive a bijection showing $F=G$.
The Garsia-Milne approach is to multiply by the reciprocal power series
$K = H^{-1}$.
A bijection showing $H \cdot K=1$ yields bijections showing that
$F=F \cdot H \cdot K$ and $G=G \cdot H \cdot K$.
Now we have a chain of bijections showing
\[
F = F \cdot H \cdot K = G \cdot H \cdot K = G
\pd
\]
In each degree, we have a chain of bijections of signed sets,
beginning and ending with unsigned sets.
By subtraction, we get a matching between the terms of $F$ and $G$.

This clever way to sidestep division was a key aspect
of Garsia and Milne's work.
\begin{comment}
We celebrate it with an anagram:

\vspace{0.5cm}
\centerline{
UNINHABITED SETS PROCRASTINATE OVER IDIOT THEOREM.
}
\end{comment}

\newpage

\section{Koenig's proof} \label{app}

Here is Koenig's proof of the Cantor-Schroeder-Bernstein equivalence theorem,
reprinted from \cite{koenig},
with some trivial misprints corrected.

\vspace{0.5cm}

\centerline{{\sc S\'eance du 9 Juillet 1906.}}

\vspace{0.5cm}

\centerline{{\sc Analyse Math\'ematique} ---
\emph{Sur la th\'eorie des ensembles.}}
\centerline{Note de {\bf M. Jules Koenig},
pr\'esent\'ee par M. H. Poincar\'e.}

\vspace{0.5cm}

La nouvelle d\'emonstration du th\'eor\`eme d'\'equivalence de M. Cantor
que je veux donner dans ces lignes a,
comme je crois,
une importance assez grande,
vu la discussion actuelle sur les fondements de la logique,
de l'arithm\'etique
et de la th\'eorie des ensembles.
Je ne voulais la donner
que dans l'exposition de la
\emph{Logique synth\'etique},
que j'esp\`ere publier bient\^ot
et que j'ai d\'eja donn\'ee dans mon cours de cette ann\'ee.
Mais l'int\'er\^et qu'on prend aujourd'hui \`a ces choses
me fait publier cette Note.

La critique spirituelle et profonde de M. Poincar\'e
(voir la \emph{Revue de M\'etaphysique et de Morale}, mai 1906)
est irr\'efutable,
\`a ce que je crois,
dans ses parties n\'egatives.
Ce que nous poss\'edons jusqu'\`a pr\'esent
\'etait peut-\^etre n\'ecessaire pour le d\'eveloppement
de la nouvelle science logique;
mais certainement cela ne donne pas ce que nous cherchons:
les bases de cette nouvelle science.

Quant au th\'eor\`eme cit\'e,
\'enonc\'e pour la premi\`ere fois par M. Cantor
et d\'emontr\'e apr\`es par MM.\ Bernstein,
Schroeder et Zermelo,
il nous faudrait le mettre en \'evidence,
sans employer le concept de nombre.

De plus nous devrions \'eviter le principe d'induction compl\`ete,
pendant que,
comme M. Poincar\'e l'a remarqu\'e bien justement,
toutes les d\'emonstrations publi\'ees jusqu'ici en font emploi.
(Quant au concept de nombre,
il est bien vrai que nous devons le construire nous-m\^emes.
Il y en a bien quelque chose dans l'intuition imm\'ediate,
un \emph{fait v\'ecu} ou une \emph{\'experience};
mais ce r\'esidu est de toute n\'ecessit\'e.)

Le th\'eor\`eme d'\'equivalence est un th\'eor\`eme d'intuition.
Pour d\'emontrer cela j'emploierai la terminologie de M. Cantor;
mais en soulignant en m\^eme temps qu'une exposition plus \'etendue
et plus pr\'ecise ne pourrait plus
se servir des mots \emph{ensemble}, etc.

{\it
Soient $X,Y$ des ensembles d\'etermin\'es,
$X_1,Y_1$ des ensembles partiels de $X$ et de $Y$ respectivement.
Nous devons d\'emontrer que,
\'etant $X \sim Y_1$ et $Y \sim X_1$,
nous aurons toujours $X \sim Y$.
}

La proposition $X \sim Y_1$ signifie la supposition de la loi (I) suivante:

{\it
Un \'el\'ement quelconque $x$ de $X$ d\'etermine un et un seul \'el\'ement
$y$ de $Y$;
donc cet $y$ d\'etermine aussi le $x$ correspondant.
Mail il y a un ou plusieurs \'el\'ements de $Y$ qui ne figurent pas
dans cette loi.
}

De m\^eme la proposition $Y \sim X_1$
signifie la supposition d'une loi (II),
qu'il serait superflu de d\'etailler encore.

Prenons donc un \'el\'ement quelconque $x_1$ de $X$;
apr\`es (I),
il nous donne un \'el\'ement d\'etermin\'e
$y_1$ de
%$Y$;
$Y_1$;
cet \'el\'ement $y_1$ nous donne,
puis par la loi (II),
un \'el\'ement d\'etermin\'e $x_2$ de $X_1$, etc.
En faisant cela,
nous ne \emph{comptons} pas;
il n'y a l\`a qu'un emploi des signes $1,2,\ldots$ pour distinguer
les \'el\'ements de $X$.
Mais les concepts \emph{suivre} et \emph{suite} doive bien \^etre
accept\'es comme concepts logiques d\'efinitifs.

Ainsi la suite
\[
x_1 y_1 x_2 y_2 \ldots
\]
peut toujours \^etre continu\'ee \`a droite,
mais pas toujours \`a gauche.
Si $x_1$ est un \'el\'ement de $X_1$,
la loi (II) donne un \'el\'ement $y_0$,
qui pr\'ec\`ede imm\'ediatement
%$x$;
$x_1$;
mais si $x_1$ est un \'el\'ement de $X$,
qui ne se trouve pas dans $X_1$,
la suite ne pourra plus \^etre continu\'ee \`a gauche.

On voit donc que les cas possibles sont les suivants:

La suite commence avec un \'el\'ement de $X$.
La suite commence avec un \'el\'ement de $Y$.
La suite peut toujours \^etre prolong\'ee \`a gauche.

Les \'el\'ements $x'_1$ et $x''_1$ de $X$ nous donnent ainsi deux suites
correspondantes:
\[
\tag{1} x'_1 y'_1 x'_2 y'_2 \ldots
\]
\[
\tag{2} x''_1 y''_1 x''_2 y''_2 \ldots
.
\]

S'il y a un \'el\'ement commun dans las suites (1) et (2),
l'\'el\'ement qui le suit est d\'etermin\'e par la loi (I),
en cons\'equence il sera le m\^eme dans les
suites (1) et (2), de m\^eme le pr\'ec\'edent s'il y en a.

C'est-\`a-dire:
Un \'el\'ement quelconque de $X$ d\'etermine toujours la suite correspondante.
Il n'est pas n\'ecessaire de d\'etailler le cas sp\'ecial d'une
suite \emph{p\'eriodique}.
C'est \'evident,
qu'une suite p\'eriodique peut toujours \^etre
prolong\'ee \`a gauche.

La loi d'\'equivalence,
dont l'expression est $X \sim Y$,
se trouve d\'etermin\'ee de fait
par ces consid\'erations.

Soit $\bar{x}$ un \'el\'ement quelconque de $X$;
nous avons l'instruction
pour la formation de la suite correspondante.
Si cette suite commence avec un \'el\'ement de $X$,
ou si elle peut \^etre continu\'ee \`a gauche,
nous choisirons comme \'el\'ement correspondant \`a $\bar{x}$ dans $Y$
l'\'el\'ement qui le suit dans la suite.
Mais, si la suite commence avec un \'el\'ement de $Y$,
nous prendrons comme \'el\'ement correspondant dans $Y$ celui qui pr\'ec\`ede
$\bar{x}$ imm\'ediatement.

Ainsi l'\'equivalence $X \sim Y$ est fix\'ee.
L'intuition pure nous m\`ene \`a reconnaitre son \emph{existence}.

Il va sans dire que cette exposition a encore beaucoups d'inconv\'enients;
parce que nous n'avons pas discut\'e \`a fond les concepts logiques
qui s'y trouve.
Telle est aussi l'expression \emph{\`a droite} ou \emph{\`a gauche}.

\newpage
\bibliography{zzz}
\bibliographystyle{hplain}

\end{document}